# A MIXTURE REPRESENTATION OF $\pi$ WITH APPLICATIONS IN MARKOV CHAIN MONTE CARLO AND PERFECT SAMPLING

By James P. Hobert[1] and Christian P. Robert[2]

*University of Florida, and Université Paris Dauphine and CREST, INSEE*

Let $X = \{X_n : n = 0, 1, 2, \ldots\}$ be an irreducible, positive recurrent Markov chain with invariant probability measure $\pi$. We show that if $X$ satisfies a one-step minorization condition, then $\pi$ can be represented as an infinite mixture. The distributions in the mixture are associated with the hitting times on an accessible atom introduced via the splitting construction of Athreya and Ney [*Trans. Amer. Math. Soc.* **245** (1978) 493–501] and Nummelin [*Z. Wahrsch. Verw. Gebiete* **43** (1978) 309–318]. When the small set in the minorization condition is the entire state space, our mixture representation of $\pi$ reduces to a simple formula, first derived by Breyer and Roberts [*Methodol. Comput. Appl. Probab.* **3** (2001) 161–177] from which samples can be easily drawn. Despite the fact that the derivation of this formula involves no coupling or backward simulation arguments, the formula can be used to reconstruct perfect sampling algorithms based on coupling from the past (CFTP) such as Murdoch and Green's [*Scand. J. Statist.* **25** (1998) 483–502] Multigamma Coupler and Wilson's [*Random Structures Algorithms* **16** (2000) 85–113] Read-Once CFTP algorithm. In the general case where the state space is not necessarily 1-small, under the assumption that $X$ satisfies a geometric drift condition, our mixture representation can be used to construct an arbitrarily accurate approximation to $\pi$ from which it is straightforward to sample. One potential application of this approximation is as a starting distribution for a Markov chain Monte Carlo algorithm based on $X$.

Received October 2002; revised May 2003.

[1]Supported in part by NSF Grant DMS-00-72827 and by the Centre National de la Recherche Scientifique (CNRS) during a visit to the Centre de Recherche en Economie et Statistique (CREST) in Paris in June 2000.

[2]Supported in part by the Department of Statistics at the University of Florida during a visit in August 1999.

**1. Representing $\pi$ as a mixture.** Let $P(x, dy)$ be a Markov transition kernel on a general state space $(\mathsf{X}, \mathcal{B}(\mathsf{X}))$ and write the associated discrete time Markov chain as $X = \{X_n : n = 0, 1, 2, \dots\}$. For $t \in \mathbb{N} := \{1, 2, 3, \dots\}$, let $P^t(x, dy)$ denote the $t$-step Markov transition kernel corresponding to $P$. Then for $n \in \mathbb{N}$, $x \in \mathsf{X}$ and a measurable set $A$, $P^t(x, A) = \Pr(X_{t+n} \in A | X_n = x)$. Throughout the paper we assume that $X$ is $\pi$-irreducible and positive Harris recurrent where $\pi$ is the invariant probability measure; for definitions see Meyn and Tweedie [(1993), Part I]. For an arbitrary measure $\mu$ and function $g$, we use the usual notation $P^t(\mu, A) = \int_{\mathsf{X}} P^t(x, A) \mu(dx)$ and $\mu(g) = \int_{\mathsf{X}} g(x) \mu(dx)$.

The assumptions we have made guarantee the existence of an $m \in \mathbb{N}$, a probability measure $\nu$ on $\mathcal{B}(\mathsf{X})$, a *small set* $C$ with $\pi(C) > 0$ and an $\varepsilon > 0$ such that for any $x \in C$,

$$P^m(x, A) \geq \varepsilon \nu(A) \qquad \forall A \in \mathcal{B}(\mathsf{X}).$$

For ease of exposition, we consider only the *strongly aperiodic* case in which $m = 1$; that is, we assume $X$ satisfies a *one-step minorization condition*

$$(1) \qquad P(x, \cdot) \geq \varepsilon \nu(\cdot) \qquad \forall x \in C.$$

A minorization condition allows for the celebrated *splitting construction* of Athreya and Ney (1978) and Nummelin (1978, 1984). To be specific, if $x \in C$, we can use (1) to write $P(x, \cdot)$ as a two-component mixture

$$(2) \qquad P(x, \cdot) = \varepsilon \nu(\cdot) + (1 - \varepsilon) R(x, \cdot),$$

where $R(x, \cdot) = (1 - \varepsilon)^{-1}(P(x, \cdot) - \varepsilon \nu(\cdot))$ is called the *residual measure*; define $R(x, \cdot)$ to be 0 if $\varepsilon = 1$.

If $X_n \in C$, then (2) can be used to generate $X_{n+1}$ sequentially as follows. Draw $\delta_n \sim \text{Ber}(\varepsilon)$. If $\delta_n = 1$, then take $X_{n+1}$ from $\nu(\cdot)$, else take $X_{n+1}$ from $R(X_n, \cdot)$. This can be formalized by introducing the *split* chain, $X' = \{(X_n, \delta_n) : n = 0, 1, \dots\}$, which lives on the space $\mathsf{X} \times \{0, 1\}$ and has Markov transition kernel given by

$$P'[(x, 0), A \times \{\delta\}] = \begin{cases} [\varepsilon \delta + (1 - \varepsilon)(1 - \delta)] P(x, A), & \text{if } x \notin C, \\ [\varepsilon \delta + (1 - \varepsilon)(1 - \delta)] R(x, A), & \text{if } x \in C, \end{cases}$$

and

$$P'[(x, 1), A \times \{\delta\}] = \begin{cases} [\varepsilon \delta + (1 - \varepsilon)(1 - \delta)] P(x, A), & \text{if } x \notin C, \\ [\varepsilon \delta + (1 - \varepsilon)(1 - \delta)] \nu(A), & \text{if } x \in C, \end{cases}$$

where $\delta \in \{0, 1\}$. It is clear that, marginally, the sequence $\{X_n : n = 0, 1, \dots\}$ from the split chain is equivalent (in distribution) to $X$. Also, the measure $\pi'$ defined by

$$\pi'(A \times \{\delta\}) = \pi(A)[\varepsilon \delta + (1 - \varepsilon)(1 - \delta)]$$

is invariant for $X'$. The key to our argument is that the set $\alpha := C \times \{1\}$ is an *accessible atom* [Meyn and Tweedie (1993), page 100] and the (random) times at which $X'$ enters $\alpha$ are *regeneration times*. Define $\tau_\alpha$ to be the first return time to $\alpha$; that is,

$$\tau_\alpha = \min\{n \geq 1 : (X_n, \delta_n) \in \alpha\}.$$

Also, let $\Pr_\alpha(\cdot)$ and $E_\alpha(\cdot)$ denote probability and expectation conditional on $(X_0, \delta_0) \in \alpha$; that is, conditional on $X_1 \sim \nu(\cdot)$. It follows from Kac's theorem [Meyn and Tweedie (1993), Theorem 10.2.2] that $E_\alpha(\tau_\alpha) < \infty$. Hence, we may define a nonnegative, nonincreasing sequence $\{p_t\}_{t=1}^\infty$ that sums to one by putting

$$(3) \qquad p_t = \frac{\Pr_\alpha(\tau_\alpha \geq t)}{E_\alpha(\tau_\alpha)}.$$

Now, for any $t \in \mathbb{N}$ and any measurable set $A$, define

$$(4) \qquad Q_t(A) = \Pr_\alpha(X_t \in A | \tau_\alpha \geq t);$$

that is, $Q_t$ is the conditional distribution of $X_t$ given that $(X_0, \delta_0) \in \alpha$ and that there are no regenerations in the split chain before time $t$. We now state the first of our two main results.

THEOREM 1. *Let $P$ be a Markov transition kernel on a general state space $(\mathsf{X}, \mathcal{B}(\mathsf{X}))$. Assume that the associated Markov chain, $X$, is $\pi$-irreducible and positive Harris recurrent where $\pi$ is the invariant probability measure. Assume further that the minorization condition* (1) *holds. Then for any $A \in \mathcal{B}(\mathsf{X})$, we have*

$$(5) \qquad \pi(A) = \sum_{t=1}^\infty Q_t(A) p_t,$$

*where $p_t$ and $Q_t(\cdot)$ are defined at* (3) *and* (4), *respectively.*

PROOF. Applying Meyn and Tweedie's (1993) Theorem 10.2.1 to $X'$ and using the fact that $\pi'(A \times \{0, 1\}) = \pi(A)$, we have

$$\pi(A) = \frac{1}{E_\alpha(\tau_\alpha)} \sum_{t=1}^\infty \Pr_\alpha(X_t \in A, \tau_\alpha \geq t) = \sum_{t=1}^\infty \Pr_\alpha(X_t \in A | \tau_\alpha \geq t) p_t. \qquad \square$$

Representation (5) is appealing from a simulation point of view because it reveals the potential for drawing from $\pi$ by randomly drawing an element from the set $\{Q_1, Q_2, Q_3, \ldots\}$ according to the probabilities $p_1, p_2, p_3, \ldots$ and then making an independent random draw from the chosen $Q_t$. This idea is closely related to *perfect sampling* [Fill (1998) and Propp and Wilson (1996)], which is a simulation method wherein a Markov chain with stationary distribution $\pi$ is used to produce independent and identically distributed (i.i.d.) samples from $\pi$. In fact, it is shown in Section 2 that when $C = \mathsf{X}$, there are direct connections between (5) and perfect sampling. We show in Section 3 that if $X$ also satisfies a drift condition, our mixture representation can be used to construct an arbitrarily accurate approximation of $\pi$ from which it is easy to sample. Finally, in Section 4 we explain how this approximation to $\pi$ provides a new method of dealing with the *burn-in* problem in Markov chain Monte Carlo (MCMC).

**2. The case where $C = \mathsf{X}$: perfect sampling.** Consider the special case in which $C = \mathsf{X}$, which of course implies that $X$ is 1-*uniformly ergodic*. In this case, $X_n \in C$ for all $n$ and, hence, $\tau_\alpha \sim \text{Geo}(\varepsilon)$; that is, $\Pr_\alpha(\tau_\alpha = t) = \varepsilon(1-\varepsilon)^{t-1}$ for $t \in \mathbb{N}$. Plugging into (3) yields $p_t = \varepsilon(1-\varepsilon)^{t-1}$ so the $p_t$s are also geometric probabilities. The distribution $Q_t$ is also quite simple when $C = \mathsf{X}$. Indeed, if there were no regenerations in the split chain before time $t$, this means that, after having drawn $X_1 \sim \nu$, the residual measure, $R$, was applied $t-1$ consecutive times to get $X_t$. We state this as a corollary.

COROLLARY 1.  *Under the assumptions of Theorem* 1, *if* $C = \mathsf{X}$, *then*

$$(6) \qquad \pi(A) = \sum_{t=1}^{\infty} \varepsilon(1-\varepsilon)^{t-1} R^{t-1}(\nu, A),$$

*where* $R^0(\nu, A)$ *is defined as* $\nu(A)$.

It is possible to derive (6) directly without using (5). Indeed, consider a Markov transition kernel, $M$, on the general space $(\mathsf{Z}, \mathcal{B}(\mathsf{Z}))$ that takes the form $M(z, \cdot) = \omega \mu(\cdot) + (1-\omega)K(z, \cdot)$, where $K(z, \cdot)$ is a $\psi$-irreducible Markov transition kernel on the same space, $\mu(\cdot)$ is a probability measure on $\mathcal{B}(\mathsf{Z})$ and $\omega \in (0,1)$. Theorem 2 of Breyer and Roberts (2001) shows that $\sum_{j=1}^{\infty} \omega(1-\omega)^{j-1} K^{j-1}(\mu, \cdot)$ is an invariant probability measure for $M$. Now note that when (1) holds with $C = \mathsf{X}$, then $P(x, \cdot) = \varepsilon \nu(\cdot) + (1-\varepsilon)R(x, \cdot)$ for all $x \in \mathsf{X}$ and it follows that $\pi$ can be written in the form (6). (It is interesting to note that $M$ is positive recurrent even if $K$ is badly behaved, e.g., transient.)

Corollary 1 immediately yields the following algorithm for sampling from $\pi$:

ALGORITHM I.

1. Simulate $x_1 \sim \nu(\cdot)$ and, independently, $t \sim \text{Geo}(\varepsilon)$.
2. If $t = 1$, take $x_1$, else simulate the transition $x_{n+1} \sim R(x_n, \cdot)$ for $n = 1, \ldots, t-1$ and take $x_t$.

Algorithm I is exactly the Multigamma Coupler of Murdoch and Green [(1998), page 486] which is a perfect sampling algorithm based on *coupling from the past* (CFTP) [Propp and Wilson (1996)]. Note that we have used our mixture representation of $\pi$ to derive this algorithm with no appeal to coupling or backward simulation.

Breyer and Roberts (2001) show, in the context of their catalytic perfect simulation algorithm, that Corollary 1 can also be used to derive Wilson's (2000) Read-Once CFTP algorithm. We now give a slightly different and more detailed description of this connection which culminates in a statement of the algorithm that Wilson described on page 93 of his paper [Wilson (2000)]. We begin with a Markov transition kernel, $S$, on the general state space $(\mathsf{X}, \mathcal{B}(\mathsf{X}))$, such that the associated Markov chain, $Y = \{Y_0, Y_1, \ldots\}$, is $\pi$-irreducible and positive Harris recurrent where $\pi$ is the invariant probability measure. Let $g : \mathsf{X} \times (0,1) \to \mathsf{X}$ be a function such that if $U \sim \mathrm{Uni}(0,1)$, then for any $x \in \mathsf{X}$ and any measurable $A$,

$$\Pr[g(x,U) \in A] = S(x,A).$$

Now for $n \in \mathbb{N}$ let $G_n : \mathsf{X} \times (0,1)^n \to \mathsf{X}$ be defined through compositions of $g$ as follows:

$$G_n(x, u_1, \ldots, u_n) = g(g(\cdots g(x, u_1) \cdots, u_{n-1}), u_n).$$

For example, $G_3(x, u_1, u_2, u_3) = g(g(g(x, u_1), u_2), u_3)$. Clearly, if $U_1, \ldots, U_n$ are i.i.d. $\mathrm{Uni}(0,1)$, then $G_n(x, U_1, \ldots, U_n)$ has distribution $S^n(x, \cdot)$. If $G_n(x, u_1, \ldots, u_n)$ is constant in $x$ for some fixed $(u_1, \ldots, u_n)$, we call $G_n(x, u_1, \ldots, u_n)$ a *coalescent* composite map.

REMARK 1. Wilson's setup is actually a bit more abstract than ours. First, he does not assume as much as we do about the structure of the "random function" $g$. Second, Wilson assumes that the user possesses an efficient, but *imperfect* method for checking whether $G_n(x, u_1, \ldots, u_n)$ is coalescent. This method will never incorrectly conclude coalescence, but may miss the fact that a particular $G_n(x, u_1, \ldots, u_n)$ is coalescent. When this imperfect method concludes that $G_n(x, u_1, \ldots, u_n)$ is constant in $x$ for some fixed $(u_1, \ldots, u_n)$, then $G_n(x, u_1, \ldots, u_n)$ is called "officially coalescent."

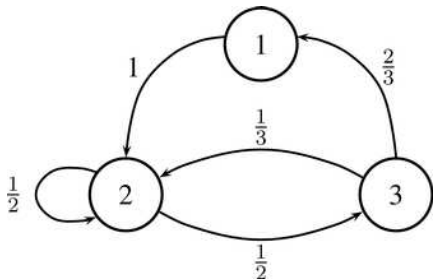

FIG. 1. *Probability transition diagram for a Markov chain on* $\mathsf{X} = \{1, 2, 3\}$.

Suppose that $k \in \mathbb{N}$ is such that when $U_1, \ldots, U_k$ are i.i.d. Uni$(0,1)$, we have

$$\text{(7)} \qquad \Pr[G_k(x, U_1, \ldots, U_k) \text{ is coalescent}] = \varepsilon > 0.$$

As an example, consider the Markov chain on $\mathsf{X} = \{1, 2, 3\}$ whose evolution is described by the probability transition diagram given in Figure 1. For this chain, we could take $g$ as follows:

$$g(1, u) = 2, \qquad g(2, u) = 2I_{(0,1/2)}(u) + 3I_{[1/2,1)}(u)$$

and

$$g(3, u) = I_{(0,2/3)}(u) + 2I_{[2/3,1)}(u).$$

It is not difficult to see that $\Pr[G_1(x, U_1) \text{ is coalescent}] = 0$, but $\Pr[G_2(x, U_1, U_2) \text{ is coalescent}] = 1/4$.

Now if we set $P = S^k$, then for any measurable set $A$,

$$\text{(8)} \qquad \begin{aligned} P(x, A) &= \Pr[G_k(x, U_1, \ldots, U_k) \in A] \\ &= \varepsilon \Pr[G_k(x, U_1, \ldots, U_k) \in A | G_k \text{ is coalescent}] \\ &\quad + (1 - \varepsilon) \Pr[G_k(x, U_1, \ldots, U_k) \in A | G_k \text{ is not coalescent}] \\ &= \varepsilon \nu(A) + (1 - \varepsilon) R(x, A), \end{aligned}$$

where we have defined $\nu(\cdot)$ and $R(x, \cdot)$ in an obvious way. Drawing from $\nu(\cdot)$ is quite simple—just simulate i.i.d. copies of $(U_1, \ldots, U_k)$ until the observed value of $G_k(x, U_1, \ldots, U_k)$ is coalescent. Drawing from $R(x, \cdot)$ can be done similarly by waiting for the first noncoalescent value of $G_k(x, U_1, \ldots, U_k)$. Therefore, if $\varepsilon$ is known, Algorithm I can be applied to make draws from $\pi$.

The beauty of (8), however, is that it can be used to simulate from $\pi$ even when $\varepsilon$ is unknown! Indeed, assume that (7) holds, but that the exact value of $\varepsilon$ is unknown. Note that application of Algorithm I does not require knowledge of $\varepsilon$, it only requires the ability to simulate from the Geo$(\varepsilon)$ distribution. But simulating $Z \sim \text{Geo}(\varepsilon)$ is easy—just generate i.i.d. copies of $(U_1, \ldots, U_k)$ and let $Z$ be the number of trials required until the first coalescent $G_k(x, U_1, \ldots, U_k)$ is observed. Moreover, the byproducts of simulating $Z$ are $Z - 1$ i.i.d. noncoalescent $G_k(x, U_1, \ldots, U_k)$s and an independent coalescent $G_k(x, U_1, \ldots, U_k)$. All of this is formalized in the following algorithm, which describes how to use (8) to simulate from $\pi$.

ALGORITHM II.

1. Simulate a coalescent $G_k(x, U_1, \ldots, U_k)$. Call its value $s$.
2. Draw independent copies of $(U_1, \ldots, U_k)$ (where the components are $k$ i.i.d. uniforms) until the observed value of $G_k(x, U_1, \ldots, U_k)$ is coalescent. Let $t$ denote the number of trials required and write the $t$ observed values of $(U_1, \ldots, U_k)$ as $(u_{i,1}, \ldots, u_{i,k})$, $i = 1, \ldots, t$.
3. Take

$$G_k(\cdots G_k(s, u_{t-1,1}, \ldots, u_{t-1,k}) \cdots, u_{1,1}, \ldots, u_{1,k}).$$

Algorithm II is exactly Wilson's (2000) Read-Once CFTP algorithm. It lends itself to iteration. Indeed, $G_k(x, u_{t,1}, \ldots, u_{t,k})$ is coalescent and is *not used* at step 3. Thus, it can be used at step 1 of the next iteration of the algorithm.

We end this section with an interesting interpretation of Corollary 1. Because $C = \mathsf{X}$, we can assume that all transitions of $X$ are made using (2). Now, each time a regeneration occurs, that is, each time a draw is made from $\nu(\cdot)$, we have to wait a $\text{Geo}(\varepsilon)$ number of iterations before the next draw from $\nu$. And, of course, the residual measure, $R$, is used in between. Thus, what Algorithm I is actually doing is returning the states immediately prior to the draws from $\nu$. Hence, an intuitive way to state Corollary 1 is as follows: The states of the Markov chain immediately prior to regenerations have distribution $\pi$. Wilson (2000) attempts to connect this to the PASTA (Poisson Arrivals See Time Averages) phenomenon from the continuous time literature.

**3. The case where $C \neq \mathsf{X}$: approximating $\pi$.** While things are more difficult when $C \neq \mathsf{X}$, it is still possible to make draws from the distribution $Q_t$ using a simple accept-reject algorithm. All that is required is the ability to simulate the split chain. Note that $Q_1(\cdot) = \nu(\cdot)$, so in the algorithm, it is assumed that $t \geq 2$.

ALGORITHM III.

1. Take $(x_0, \delta_0) \in \alpha$ and simulate the split chain for $t$ iterations.
2. If, for each $i = 1, 2, \ldots, t-1$, $(x_i, \delta_i) \notin \alpha$, then take $x_t$; otherwise, repeat.

Jones and Hobert (2001) provide some practical advice on simulating the split chain.

Let $T$ denote a discrete random variable with support $\mathbb{N}$ such that $\Pr(T = t) = p_t$. The ability to randomly draw an element from the set $\{Q_1, Q_2, Q_3, \ldots\}$ according to the probabilities $p_1, p_2, p_3, \ldots$ is tantamount to being able to simulate $T$. While $T \sim \text{Geo}(\varepsilon)$ when $C = \mathsf{X}$, its distribution is quite complicated when $C \neq \mathsf{X}$. On the other hand, making i.i.d. draws from the distribution of $\tau_\alpha$ is straightforward—just take $X_1 \sim \nu(\cdot)$, run the split chain, and count the number of iterations until the first regeneration. Unfortunately, despite the simple relationship between their mass functions, it is not clear how to use i.i.d. draws from the distribution of $\tau_\alpha$ to get i.i.d. draws from the distribution of $T$. Hence, we focus on using (5) to construct an approximation to $\pi$.

Let $\{\tilde{p}_t\}_{t=1}^\infty$ denote another nonnegative sequence that sums to one and let $\tilde{T}$ denote the corresponding discrete random variable, that is, $\Pr(\tilde{T} = t) = \tilde{p}_t$. Consider an approximation of $\pi$ given by

$$\tilde{\pi}(\cdot) = \sum_{t=1}^\infty Q_t(\cdot)\tilde{p}_t. \tag{9}$$

Note that

$$\|\pi(\cdot) - \tilde{\pi}(\cdot)\| = \left\|\sum_{t=1}^\infty Q_t(\cdot)\, p_t - \sum_{t=1}^\infty Q_t(\cdot)\, \tilde{p}_t\right\| \le \sum_{t=1}^\infty |p_t - \tilde{p}_t|.$$

Thus, the total variation distance between the distributions $\pi$ and $\tilde{\pi}$ is bounded above by twice the total variation distance between the distributions of $T$ and $\tilde{T}$.

We now show that, under an additional assumption on $X$, for any given $\gamma > 0$, it is possible to construct a sequence $\{\tilde{p}_t\}_{t=1}^\infty$ such that $\sum_{t=1}^\infty |p_t - \tilde{p}_t| \le \gamma$ *and* such that making i.i.d. draws from the distribution of $\tilde{T}$ is straightforward. The assumption is that the Markov chain $X$ satisfies a *geometric drift condition*, that is, for some function $V : \mathsf{X} \to [1, \infty)$, some $\lambda < 1$ and some $b < \infty$, we have

$$PV(x) \le \lambda V(x) + b I_C(x) \qquad \forall\, x \in \mathsf{X}, \tag{10}$$

where $PV(x) := \int_\mathsf{X} V(y)\, P(x, dy)$. It is well known that (10) combined with the smallness of the set $C$ implies that the Markov chain $X$ is geometrically ergodic [Meyn and Tweedie (1993), Chapter 15]. We will need the following result, which is a combination of Theorems 2.3 and 4.1 in Roberts and Tweedie (1999).

THEOREM 2 [Roberts and Tweedie (1999)]. *Let $P$ be a Markov transition kernel on a general state space $(\mathsf{X}, \mathcal{B}(\mathsf{X}))$. Assume that the associated Markov chain, $X$, is $\pi$-irreducible and positive Harris recurrent where $\pi$ is the invariant probability measure. Assume further that the minorization condition* (1) *and the drift condition* (10) *both hold. Define $d = \sup_{x \in C} V(x)$, $A = \sup_{x \in C} PV(x)$ and $J = (A - \varepsilon)/\lambda$. Then the generating function $E_\alpha(\beta^{\tau_\alpha})$ converges for $\beta \in (1, \beta^*)$ where*:

1. *If $J < 1$, then $\beta^* = \lambda^{-1}$.*

2. *If $J \geq 1$, then*

$$\beta^* = \exp\left\{\frac{\log \lambda \log(1-\varepsilon)}{\log J - \log(1-\varepsilon)}\right\} \leq \lambda^{-1}.$$

*Furthermore, letting $\phi(\beta) = \log \beta / \log \lambda^{-1}$, if $\beta \in (1, \beta^*)$, then*

$$\Pr_\alpha(\tau_\alpha \geq t) \leq \beta[\nu(V)]^{\phi(\beta)}\left[\frac{1-\beta(1-\varepsilon)}{1-(1-\varepsilon)(J/(1-\varepsilon))^{\phi(\beta)}}\right]\beta^{-t}$$
(11)
$$=: g(\beta, \varepsilon, J)\beta^{-t}.$$

REMARK 2. In applications $\nu(V)$ may be difficult to calculate. Fortunately, there is a simple upper bound. Indeed, an application of Lemma 1 from Hobert, Jones, Presnell and Rosenthal (2002) yields $\nu(V) \leq \pi(V)/[\varepsilon\pi(C)]$. From (10) we know that $\pi(V)/\pi(C) \leq b/(1-\lambda)$ and, hence, $\nu(V) \leq b/[\varepsilon(1-\lambda)]$.

We are now in a position to state the second of our two main results.

THEOREM 3. *Assume the hypotheses of Theorem 2. Fix $\gamma > 0$ and $\beta \in (1, \beta^*)$. Let $\tilde{T}$ be the random variable supported on $\{1, \ldots, M\}$ with probabilities*

$$\tilde{p}_t = \frac{\Pr_\alpha(\tau_\alpha \geq t)}{\sum_{s=1}^M \Pr_\alpha(\tau_\alpha \geq s)},$$

*where $M$ is any integer larger than*

$$\log\left[\frac{2\,g(\beta,\varepsilon,J)}{\gamma(\beta-1)}\right]\bigg/\log \beta.$$

*Then $\|\pi(\cdot) - \tilde{\pi}(\cdot)\| \leq \gamma$.*

PROOF. First,

$$\sum_{t=1}^\infty |p_t - \tilde{p}_t| = \sum_{t=1}^M \left|\frac{\Pr_\alpha(\tau_\alpha \geq t)}{\mathrm{E}_\alpha(\tau_\alpha)} - \frac{\Pr_\alpha(\tau_\alpha \geq t)}{\sum_{s=1}^M \Pr_\alpha(\tau_\alpha \geq s)}\right| + \sum_{t=M+1}^\infty \frac{\Pr_\alpha(\tau_\alpha \geq t)}{\mathrm{E}_\alpha(\tau_\alpha)}$$

$$= 2\sum_{t=M+1}^\infty \frac{\Pr_\alpha(\tau_\alpha \geq t)}{\mathrm{E}_\alpha(\tau_\alpha)}.$$

Thus, since $\mathrm{E}_\alpha(\tau_\alpha) \geq 1$, it suffices to show that $\sum_{t=M+1}^\infty \Pr_\alpha(\tau_\alpha \geq t) \leq \gamma/2$. Using (11) from Theorem 2, we have

$$\sum_{t=M+1}^\infty \Pr_\alpha(\tau_\alpha \geq t) \leq g(\beta,\varepsilon,J)\sum_{t=M+1}^\infty \beta^{-t} = \frac{g(\beta,\varepsilon,J)}{\beta-1}\beta^{-M},$$

and the result follows from the assumption on $M$. $\square$

Of course, Theorem 3 is useful from a practical standpoint only if it is possible to sample from the distribution of $\tilde{T}$. To this end, consider the random vector $(V, W)$, where $V$ and $W$ are independent, $V$ is uniform on $\{1, \ldots, M\}$, and $W$ is equal in distribution to $\tau_\alpha$ when $(X_0, \delta_0) \in \alpha$. Note that, for any $t \in \{1, \ldots, M\}$, we have

$$\Pr(V = t | W \geq V) = \frac{\Pr(W \geq V | V = t) M^{-1}}{\sum_{i=1}^M \Pr(W \geq V | V = i) M^{-1}} = \frac{\Pr_\alpha(\tau_\alpha \geq t)}{\sum_{i=1}^M \Pr_\alpha(\tau_\alpha \geq i)} = \tilde{p}_t.$$

Hence, the following algorithm can be used to sample from the distribution of $\tilde{T}$.

ALGORITHM IV.

1. Draw $v \sim \text{Uni}\{1, \ldots, M\}$ and, independently, draw $w$ from the distribution of $\tau_\alpha$ with $(x_0, \delta_0) \in \alpha$.
2. If $w \geq v$, take $v$; otherwise, repeat.

We conclude that, given any $\gamma > 0$, Algorithms III and IV can be used to make i.i.d. draws from $\tilde{\pi}$ satisfying $\|\pi(\cdot) - \tilde{\pi}(\cdot)\| \leq \gamma$. In the last section we briefly describe how our approximation may provide an alternative solution to the *burn-in* problem in MCMC.

**4. An application to burn-in.** Suppose that the Markov kernel, $P$, is the basis of an MCMC algorithm whose purpose is to explore $\pi$. Our assumptions about $P$ imply that for every initial probability measure $\mu(\cdot)$ we have

$$\|P^n(\mu, \cdot) - \pi(\cdot)\| \downarrow 0 \quad \text{as } n \to \infty.$$

Typically, the MCMC user has no particular starting distribution in mind. Indeed, $\mu(\cdot)$ is usually taken to be a point mass at some point from which it is convenient to start the simulation. An important problem in the implementation of MCMC algorithms is *burn-in* (time), which is formally described as follows. Given $\mu(\cdot)$ and $\gamma > 0$, we want to find a value $n^*$ such that

(12) $$\|P^{n^*}(\mu, \cdot) - \pi(\cdot)\| < \gamma.$$

If (12) holds, then the marginal distribution of $X_n$ (conditional on $X_0 \sim \mu$) is within $\gamma$ of $\pi$ for all $n \geq n^*$. Hence, $n^*$ may be regarded as a reasonable time to start sampling the Markov chain.

Several authors have recently shown that drift and minorization conditions on the Markov chain can be used to derive computable upper bounds on $\|P^n(\mu, \cdot) - \pi(\cdot)\|$ that decrease geometrically fast in $n$ [Douc, Moulines and Rosenthal (2004), Meyn and Tweedie (1994), Roberts and Tweddie

A MIXTURE REPRESENTATION OF $\pi$ 11(1999) and Rosenthal (1995)]. These upper bounds can be used to find an $n^*$ that satisfies (12). Unfortunately, when this strategy is used for nontoy MCMC algorithms, it is not unusual for the resulting $n^*$s to be too large to be of any practical value [see, e.g., Jones and Hobert (2004)].

Alternatively, a seemingly unnatural way to phrase the burn-in question is as follows. Can we find a starting distribution, $\mu(\cdot)$, that is within $\gamma$ of $\pi$ in total variation? If so, we could start sampling the chain immediately. This seems unnatural because the stationary distribution of an MCMC algorithm is typically intractable and, hence, not easily approximated. Nevertheless, the results in the previous section show that we can actually construct such a starting distribution. An alternative method of dealing with the burn-in problem is to start the chain by drawing $X_0 \sim \tilde{\pi}$ and using all the samples right from the start.

Of course, $\tilde{\pi}$ would normally be constructed using the same drift and minorization conditions that are used to construct the upper bounds mentioned above. One might suspect that in situations where the $n^*$s calculated using the upper bounds are too large, simulating from $\tilde{\pi}$ might be extremely inefficient, perhaps to the point where it is not practical. On the other hand, $\tilde{\pi}$ was derived without using several inequalities that are required in deriving the upper bounds. For example, we did not use the *coupling inequality*, nor did we have to worry about constructing a *bivariate* drift condition using the drift on the original chain [see Roberts and Tweedie (1999), Theorem 5.2].

**Acknowledgments.** The authors are grateful to Laird Breyer, Galin Jones, Eric Moulines, Gareth Roberts and Richard Tweedie for helpful discussions, and to an anonymous referee for constructive comments and suggestions.

Department of Statistics  
University of Florida  
Gainsville, Florida 32611  
USA  
e-mail: jhobert@stat.ufl.edu

Université paris Dauphine  
and CREST, INSEE  
Paris  
France  
e-mail: xian@ceremade.dauphine.fr